\magnification=\magstep1
\input amssym.def
\def \PP {{\Bbb P}}
\def \FF {{\Bbb F}}
\def \ZZ {{\Bbb Z}}
\def \Spec {{\rm Spec}\,}
\def \O {{\cal O}}
\def \Ind {{\rm Ind}}
\def \Gal {{\rm Gal}}
\def \Aut {{\rm Aut}}

\def \PGL {{\rm PGL}}
\def \SL {{\rm SL}}

\def \GL {{\rm GL}}

\def \surj {{\to\!\!\!\!\!\to}} 
\def \qed {\hfill\lower0.9pt\vbox{\hrule \hbox{\vrule height 0.2 cm    
  \hskip 0.2 cm \vrule height 0.2 cm}\hrule}}
\font\small=cmr8
\font\smit=cmti8
\font\smtt=cmtt8

\centerline{\bf Oort groups and lifting problems}

\medskip

\centerline{T.~Chinburg, R.~Guralnick, D.~Harbater\footnote{$^*$}{\small The authors were respectively supported in part by NSF Grants DMS-0500106, DMS-0653873, and DMS-0500118.
{\parindent=0pt\item{}  
{\smit 2000 Mathematics Subject Classification}.  Primary 12F10, 14H37, 20B25; Secondary 13B05, 14D15, 14H30.  
\item{}  \baselineskip=10pt
{\smit Key words and phrases}: curves, automorphisms, Galois groups, characteristic {\smit p}, lifting, Oort Conjecture.}}}

\medskip

{\baselineskip=10pt\narrower\noindent{\bf Abstract.} Let $k$ be an algebraically closed field of positive characteristic $p$. We consider which finite groups $G$
have the property that every faithful action of $G$ on a connected smooth projective curve over $k$ lifts to characteristic zero.  Oort conjectured that cyclic groups have this property.  We show that if a cyclic-by-$p$ group $G$ has this property, then $G$ must be either cyclic or dihedral, with the exception of $A_4$ in characteristic $2$.  This proves one direction of a strong form of the Oort Conjecture. \par}

\medskip  

\baselineskip=12pt

\noindent{\bf \S1. Introduction.}
\medskip
The motivation for this paper is the
following conjecture made by Oort in [Oo, I.7]:
\medskip
\noindent{\bf Conjecture 1.1} (Oort Conjecture) {\sl Every faithful action
of a cyclic group on a connected smooth projective curve $Y$ over an algebraically
closed field $k$ of positive characteristic $p$ lifts to characteristic $0$.}
\medskip
\noindent With $k$ as above, we will call a finite group $G$ an {\it Oort group} for $k$ if every faithful action of $G$ on a smooth connected projective curve $Y$
over $k$ lifts to characteristic $0$.  By such a lifting we mean an action of $G$ on a smooth projective
curve ${\cal Y}$ over a complete discrete valuation ring $R$ 
of characteristic $0$ and residue field $k$ together a $G$-equivariant
isomorphism between $Y$ and the special fibre ${\cal Y}
\times_R k$.     
Thus Oort's Conjecture is that cyclic groups are Oort groups.  The object of this paper is to make  a precise prediction about which $G$ are Oort groups and to prove one direction of this prediction, namely that all Oort groups are on the list we predict.  

Grothendieck's study of the tame fundamental group of curves in characteristic $p$ [Gr, Exp.~XIII, \S2] relies on the fact that tamely ramified covers can be lifted to characteristic $0$. Oort groups over $k$ can equivalently be characterized as groups $G$
such that every connected $G$-Galois cover of $k$-curves lifts to
characteristic $0$ (see  \S2).   This fact and Grothendieck's result 
imply that prime-to-$p$ groups are Oort groups for $k$. 

It was proved by Oort, Sekiguchi and Suwa in [OSS] (resp.\ by Green and Matignon in
[GM]) that a cyclic group $G$ is an Oort group if the order of $G$
is exactly divisible by  $p$ (resp.\ by $p^2$). The dihedral group of order $2p$ is an Oort group for
all $k$ of characteristic $p$, by a result shown in [Pa] for $p=2$ and in [BW] for odd $p$ (see Example~2.12(c,f) below).  By another result stated in [BW], the alternating group $A_4$ is an Oort group in characteristic $2$ (see Example~2.12(g)).  All of the above groups are cyclic-by-$p$ (i.e.\ extensions of a prime-to-$p$ cyclic group by a $p$-group), which is the form of an inertia group associated to a cover of $k$-curves.  

The above results suggest the following strengthening of the Oort Conjecture concerning cyclic groups:

\medskip

\noindent{\bf Conjecture 1.2.} (Strong Oort Conjecture) {\it If $k$ is an algebraically closed field of characteristic $p$, and if $G$ is a cyclic-by-$p$ group, then $G$ is an Oort group for $k$ if and only if $G$ is either a cyclic group, or a dihedral group of order $2p^n$ for some $n$, or (if $p=2$) $G$ is the alternating group $A_4$.}

\medskip

By Corollary~2.8 below, an arbitrary finite group $G$ is an Oort group for $k$ if and only if every cyclic-by-$p$ subgroup of $G$ is.  So Conjecture~1.2 would also determine precisely which finite groups are Oort groups, viz.\ those whose cyclic-by-$p$ subgroups are of the above form.  In [CGH2] we give a detailed description of this class of groups.

In this paper, we show the forward direction of Conjecture~1.2: If a cyclic-by-$p$ group $G$ is an Oort group for an algebraically closed field $k$ of characteristic $p$, then it must be of the asserted form.  This is shown in odd characteristic in Corollary~3.4, and in characteristic $2$ in Theorem~4.5.  

We also consider a local version of the above problem, in which actions of $G$ on $\Spec k[[x]]$ are considered, along with the corresponding notion of a {\it local Oort group} (see Section~2 below).   This notion is in fact closer to the focus of study in [OSS], [GM], [Pa] and [BW].  In this paper we also prove results that are local analogs of our global results; see Theorem~3.3 and Theorem~4.4.  The local result in odd characteristic is the natural analog of the global version.  In characteristic $2$ our local result is somewhat more complicated. 
We will prove in [CGH1] a stronger local result concerning a lifting obstruction defined by
Bertin in [B], and we also take up the question of when some faithful local $G$-action lifts to characteristic $0$.

\medskip

{\it Notation and terminology}:  In this paper, 
$k$ denotes an algebraically closed field of characteristic $p>0$.  
A {\it curve} $X$ over a field $F$ is a normal scheme of finite type over $F$ such
that ${\rm dim}(O_{X,x}) = 1$ for all closed points $x$ of $X$.  If
$R$ is a Dedekind ring, a {\it curve} ${\cal X}$ over $R$ is a normal scheme
together with a separated, flat morphism $X \to {\rm Spec}(R)$ of finite type whose fibres are curves.  

Suppose $G$
is a finite group, $B$ is a field or a Dedekind ring, and $V$ is a connected curve over $B$.  A $G$-{\it Galois cover} over $V$ consists of a faithful action of $G$ on a curve $U$ over $B$ and isomorphism over $B$ of $V$ with the quotient curve $U/G$.  We do not require $U$
to be connected.  The resulting finite morphism $U \to U/G = V$ is $G$-equivariant when we let $G$ act trivially on $V$.  If $H$ is a subgroup of $G$, and $U' \to V$ is an $H$-Galois cover of curves over $B$, then $\Ind_H^G U' \to X$ denotes the induced $G$-Galois cover obtained by taking $(G:H)$ disjoint copies of $U'$ indexed by coset representatives of $H$ in $G$.
If $B$ is an algebraically
closed field $k$, the fact that $U$ and $V$ are normal and $V$ is connected forces $U$ to be the normalization of $V$ in 
$V \times_U {\rm Spec}(k(V)) = {\rm Spec}(D)$, where $k(V)$
is the function field of $V$ and $D$ is an \'etale $G$-algebra over $k(V)$.

Given groups $N,H$, we denote by 
$N.H$ the semi-direct product of $N$ with $H$, relative to some action of $H$ on $N$.  We denote the cyclic group of order $n$ by $C_n$ (multiplicatively) or $\ZZ/n$ (additively).  
So a cyclic-by-$p$ group is of the form $P.C_n$ for some $n$ prime to $p$, where $P$ is a $p$-group.
The dihedral group of order $2n$ (and of degree $n$) is denoted here by $D_{2n}$.  So $D_4$ denotes the Klein four group, and $D_2$ the cyclic group of order $2$.  For a prime-power $q$, $\SL(n,q)$ denotes the group $\SL_n(\FF_q)$, and similarly for $\GL$ and $\PGL$.

The {\it Frattini subgroup} of a finite group $G$ (viz.\ the intersection of the maximal subgroups of $G$) is denoted by $\Phi(G)$.  (If $G$ is a $p$-group, $\Phi(G)$ is also the subgroup of $G$ generated by $p$-th powers and commutators.)  Given subgroups $E,H$ of a group $G$, the {\it centralizer} of $H$ in $E$ is the subgroup $C_E(H) = \{e \in E\,|\, (\forall h \in H)\, eh=he\} \subset E$ and the {\it normalizer} $H$ in $E$ is the subgroup $N_E(H) = \{e \in E\,|\, eHe^{-1}=H\} \subset E$.

\bigskip

\noindent{\bf \S2.  Oort groups and local Oort groups}

\medskip

Let $X$ be a smooth complete $k$-curve, and let $R$ be a mixed characteristic complete discrete valuation ring with residue field $k$.  There is a unique continuous algebra homomorphism from the ring $W(k)$ of Witt vectors over $k$ into $R$  which induces the identity map on residue fields, and $R$ is a finite extension of $W(k)$.   By [Gr, III, Cor.~7.4], there is a smooth complete $R$-curve $\cal X$ with closed fibre isomorphic to $X$; we call this a {\it model} of $X$ over $R$.  Let $Y \to X$ be a $G$-Galois cover.  We say that the $G$-Galois cover $Y \to X$ {\it lifts to $\cal X$} if there is a  smooth complete $R$-curve ${\cal Y}$ on which 
$G$ acts and an isomorphism between ${\cal X}$ and the quotient scheme ${\cal Y}/G$ such that the resulting $G$-Galois cover ${\cal Y} \to {\cal X}$ has closed fibre  $Y \to X$ (as a $G$-Galois cover). The general fibre of ${\cal Y}$ is ${\cal Y} \times_R F$ where $F$ is the fraction field of $R$, and the geometric general fibre of 
${\cal Y}$ is ${\cal Y} \times_R \bar F$ where $\bar F$ is an algebraic
closure of $F$.  These are smooth complete curves over $F$ and $\bar F$, respectively.  

We will later need the following well-known result.  Parts (a) and (b) are special cases 
of [dJ, Proposition~4.2] (see also the proof of [Ra, Proposition~5], and its corollary).  Part (c) is then immediate from the constancy of the arithmetic genus in a connected flat family (see
[H, Chapter III, Cor.~9.10]).\medskip 

\noindent {\bf Proposition 2.1.}  {\sl With the above notations, suppose $H$ is a subgroup of $G$.  

a) The morphism 
${\cal Y} \to {\cal Y}/H$ is an $H$-Galois cover of smooth complete curves over $R$
that lifts the $H$-Galois cover of smooth complete curves $Y \to Y/H$ over $k$.

b) If $H$ is normal in $G$ then ${\cal Y}/H \to {\cal X} = {\cal Y}/G$ is a $G/H$-Galois
cover of smooth complete curves over $R$ that lifts the $G/H$-Galois
cover $Y/H \to X = Y/G$ of smooth complete curves over $k$.  

c) The genera of $Y/H$, of the general fibre of ${\cal Y}/H$, and of the geometric general fibre of ${\cal Y}/H$, are equal.}
\medskip

We say that a $G$-Galois cover of smooth complete $k$-curves $Y \to X$ {\it lifts to characteristic $0$} if it lifts to a model $\cal X$ of $X$ over some discrete valuation ring $R$ as above.  If $\xi$ is a point of $X$, then we say that $Y \to X$ {\it lifts locally near $\xi$} if for some $R$ and $\cal X$ as above, and for every point $\eta$ of $Y$ over $\xi$, there is an $I$-Galois cover $\hat {\cal Y}_\eta \to \hat {\cal X}_\xi := \Spec \hat \O_{{\cal X},\xi}$ whose closed fibre is isomorphic to the pullback of $Y \to X$ to $\hat X_\xi := \Spec \hat \O_{X,\xi}$ as an $I$-Galois cover; here $I$ is the inertia group of $Y \to X$ at $\eta$.  This property holds trivially if $\xi$ is not a branch point.  

\medskip

\noindent{\bf Theorem 2.2.}
{\sl Let $X$ be a smooth complete $k$-curve, and let
$Y \to X$ be a $G$-Galois cover.  Then the following are equivalent:

i) $Y \to X$ lifts to characteristic $0$.

ii) For every mixed characteristic complete discrete valuation ring $R$ with residue field $k$ and every model $\cal X$ of $X$ over $R$, there is a complete discrete valuation ring  $R'$ which is a finite extension of $R$ such that $Y \to X$ lifts to the induced model ${\cal X}'$ of $X$ over $R'$.

iii) $Y \to X$ lifts locally near each branch point.}

\medskip

\noindent{\it Proof.}  The implications (ii) $\Rightarrow$ (i) $\Rightarrow$ (iii) are trivial, so it suffices to
prove the implication (iii) $\Rightarrow$ (ii).  Let $S = \{\xi_1,\dots,\xi_r\}$ be a non-empty finite set of points containing the branch locus of $Y \to X$.  Then the cover lifts locally near each point of $S$.  The lift near $\xi_i$ is defined with respect to some model ${\cal X}_i$ of $X$ over some finite extension $R_i$ of $W(k)$.  Since ${\cal X}_i$ is smooth over $R_i$ and since the residue field $k$ of $R_i$ is algebraically closed, it follows that $\hat \O_{{\cal X}_i,\xi_i}$ is isomorphic to $R_i[[t_i]]$, where $t_i$ is a uniformizer of ${\cal X}_i$ over $R_i$, lifting a uniformizer $\bar t_i$ of $X$ at $\xi_i$.   

Let $R'$ be a complete discrete valuation ring in the algebraic closure of the fraction field of $R$ into which all of the $R_i$
embed over $R$.   Let ${\cal X}'$ be the $R'$-model of $X$ induced by $\cal X$.
Then the complete local rings of ${\cal X}'$ at the points $\xi_i$ are of the form $R'[[t_i]]$, and the local liftings on the $R_i$-curves ${\cal X}_i$ induce local liftings on ${\cal X}'$.  Inducing each of these from $I_i$ to $G$ (by taking a disjoint union of copies indexed by the cosets of $I_i$ in $G$), we obtain local (disconnected) $G$-Galois covers $\hat {\cal Y}_i$ of $\hat {\cal X}_i' :=\Spec \hat \O_{{\cal X}',\xi_i}$ for each $i$.

Let $U$ be the complement of $S$ in $X$.  So $U$ is an affine $k$-curve, say $U = \Spec A$, for some $k$-algebra $A$ of finite type.  The formal completion of ${\cal X}'$ along $U$ is given by ${\cal U}' := \Spec W(A) \otimes_{W(k)} R'$.  Let $V = \Spec B$ be the inverse image of $U$ under $Y \to X$; this is $G$-Galois and \'etale over $U$, and ${\cal V} := \Spec W(B) \otimes_{W(k)} R'$ is $G$-Galois and \'etale over ${\cal U}'$.  

For each $i$, let $\xi_i^\circ$ be the generic point of $\hat X_i := 
\Spec \hat \O_{X,\xi_i}$, and let $\hat {\cal U}_i = \Spec \hat \O_{{\cal X}',\xi_i^\circ}$. There are thus natural morphisms $\hat {\cal U}_i \to \hat {\cal X}_i'$ and $\hat {\cal U}_i \to {\cal U}'$ (and we regard $\hat {\cal U}_i$ as the ``overlap'' of $\hat {\cal X}_i'$ with ${\cal U}'$ in ${\cal X}'$).  Pulling back $\hat {\cal Y}_i \to \hat {\cal X}_i'$ via $\hat {\cal U}_i \to \hat {\cal X}_i'$ yields a $G$-Galois cover of $\hat {\cal U}_i$, and so does pulling back ${\cal V} \to {\cal U}'$ via $\hat {\cal U}_i \to {\cal U}'$.  For each of these two pullbacks, the fibre over $\xi_i^\circ$ is equipped with an isomorphism to the fibre of $Y \to X$ over $\xi_i^\circ$ (as a $G$-space).  The induced isomorphism between these fibres of the two pullbacks lifts to a unique isomorphism between these two pullbacks as $G$-Galois covers, by [Se, III, \S5, Thm.~2]. 

We now apply formal patching (e.g.\ [HS, Cor.\ to Thm.~1] or [Pr, Thm.~3.4]) to the proper $R'$-curve ${\cal X}'$ and the above data.  So there is a $G$-Galois cover ${\cal Y} \to {\cal X}'$ whose restriction to $\hat {\cal X}_i'$ is isomorphic to $\hat {\cal Y}_i$; whose restriction to ${\cal U}'$ is isomorphic to ${\cal V}$; and whose closed fibre is isomorphic to $Y \to X$.  So (ii) holds.  \qed

\medskip

\noindent{\it Remark.} A similar argument, using rigid patching, was used in the proof of [GM, III, Lifting Theorem~1.3], in the case of covers whose inertia groups are cyclic of order not divisible by $p^3$, where $p = {\rm char}\, k$.  The equivalence of (i) and (iii) was proved using deformation theory in [BM, 
Th\'eor\`eme~4.6].

\medskip

Consider a local $G$-Galois cover $\hat Y \to \hat X := \Spec k[[x]]$.
Let $R$ be a discrete valuation ring which is a finite extension of the ring of Witt vectors $W(k)$, and let $\hat {\cal X} = \Spec R[[x]]$.  We say that the given cover {\it lifts to} $\hat {\cal X}$ if there is a $G$-Galois cover $\hat {\cal Y} \to \hat {\cal X}$ whose closed fibre is $\hat Y \to \hat X$ is a $G$-Galois cover.  Similarly, we say that the $G$-Galois cover $\hat Y \to \hat X$ {\it lifts to characteristic $0$} if it lifts to $\hat {\cal X} = \Spec R[[x]]$ for some discrete valuation ring which is a finite extension of $W(k)$.  

We may identify $\hat X = \Spec k[[x]]$ with $\Spec \hat{\O}_{X,\infty}$, where $X = \PP^1_k$.  Let $G=P.C$ be a cyclic-by-$p$ group; i.e.\ a semi-direct product of a $p$-group $P$ with a cyclic group $C = C_m$ of order $m$ prime to $p$.  Recall that given any $G$-Galois cover $\hat Y \to \hat X = \Spec k[[x]]$, there is a unique $G$-Galois cover $Y \to X := \PP^1_k$ whose restriction to $\Spec \hat{\O}_{X,\infty}$ agrees with $\hat Y \to \hat X$; which is tamely ramified over $0$ with ramification index equal to $m$; and which is unramified elsewhere [Ka, Thm.~1.4.1].  Here $Y \to X$ is called the {\it Katz-Gabber} cover associated to $\hat Y \to \hat X$.  

Theorem 2.2 then has the following corollary:

\medskip

\noindent{\bf Corollary 2.3.} 
{\sl Let $G$ be a cyclic-by-$p$ group and let $\hat Y \to \hat X = \Spec k[[x]]$ be a connected $G$-Galois cover.  Let $Y \to X$ be the associated Katz-Gabber cover.  

a) The cover $Y \to X$ lifts to characteristic $0$ if and only if $\hat Y \to \hat X$ lifts.

b) Let $g$ be the genus of $Y$.  If there is no connected genus $g$ curve $Y^\circ$ over an algebraically closed field of characteristic $0$ together with a faithful action of $G$ such that $Y^\circ/G$ has genus $0$, then $\hat Y \to \hat X$ does not lift to characteristic $0$.

c)  Suppose $Y$ has genus $0$.  If there is no algebraically closed field $L$ of characteristic $0$
such that $G$ embeds into ${\rm PGL}_2(L)$ then $\hat Y \to \hat X$ does not lift to characteristic $0$.}

\medskip

\noindent{\it Proof.}  a) If $Y \to X$ lifts to characteristic $0$, then
it lifts locally by Theorem 2.2, so 
$\hat Y \to \hat X$ lifts to characteristic $0$.  

Conversely, suppose that $\hat Y \to \hat X$ lifts to characteristic $0$.  Then $Y \to X$ lifts locally near the branch point $\infty$.  But $Y \to X$ also lifts locally near the branch point $0$ since it is tamely ramified there, and tame covers lift [Gr, Exp.~XIII, \S2].  So by Theorem 2.2, $Y \to X$  lifts to characteristic $0$.  

b) If $\hat Y \to \hat X$ lifts to characteristic $0$, then so does $Y \to X$, by part (a).  Let ${\cal Y} \to {\cal X}$ be a lift to characteristic $0$, with geometric generic fibre $Y^\circ \to X^\circ$.  Since $k$ is algebraically closed and $Y$ is connected, $Y^\circ$ must be connected.  By Proposition 2.1, $X$ and $Y$ have the same genera as their generic fibres, viz.\ $0$ and
$g$ respectively.  Hence the same is true for $X^\circ$ and $Y^\circ$.
But $X^\circ = Y^\circ/G$.  This contradicts the hypothesis.

c) Let $g = 0$ in part (b), so that if $\hat Y \to \hat X$ lifts to characteristic $0$ there
is a connected genus $0$ curve $Y^\circ$ over an algebraically closed field $L$ of characteristic $0$ for which $G$ acts faithfully
on $Y^{\circ}$. This $Y^{\circ}$ must be isomorphic to $\PP^1_L$, so $G$ embeds into
${\rm Aut}_L(Y^{\circ})  = {\rm PGL}_2(L)$, which proves (c). 
\qed

\medskip

Let $G$ be a finite group.  We say that $G$ is an {\it Oort group for} $k$ if for every smooth connected complete $k$-curve $X$, every connected $G$-Galois cover $Y \to X$ lifts to characteristic $0$.  (If the field $k$ is understood, we will sometimes omit the words ``for $k$".  As F. Pop has noted, it is a very interesting question whether
the set of Oort groups for $k$ depends only on the characteristic of $k$.)   Recall that every finite group is the Galois group of some connected cover of $X$ (and moreover the absolute Galois group of the function field of $X$ is free profinite of rank ${\rm card}\, k$; cf.\ [Ha1], [Po]).  So this condition on $G$ is not vacuous.  
Note also that $Y \to X = Y/G$ lifts to characteristic $0$
if and only we may lift the action of $G$ on $Y$ to an action of $G$ on a smooth
complete curve ${\cal Y}$ over a complete discrete valuation ring $R$ of characteristic $0$ and residue
field $k$.  For if such a ${\cal Y}$ exists, the curve  ${\cal X} = {\cal Y}/G$ over $R$ will have special fibre
$({\cal Y}\times_R k)/G = Y/G = X$.  

A finite group $G$ is the Galois group of a connected cover of ${\rm Spec}(k[[x]])$ if and only if it is cyclic-by-$p$.  So if $G$ a cyclic-by-$p$ group, we will say that $G$ is a {\it local Oort group} for $k$ if every connected $G$-Galois cover of $\Spec k[[x]]$ lifts to characteristic $0$.

\medskip

\noindent{\bf Theorem 2.4.} 
{\sl Let $G$ be a finite group.  Then the following are equivalent:

i) $G$ is an Oort group for $k$.

ii) Every $G$-Galois cover of $\PP^1_k$ lifts to $\PP^1_R$, for some finite extension $R$ of $W(k)$ (depending on the cover).

iii) Every cyclic-by-$p$ subgroup of $G$ is a local Oort group for $k$.}
 
\medskip

The key step in proving this result is

\medskip

\noindent{\bf Lemma 2.5.} 
{\sl Let $G$ be a finite group, and let $I \subset G$ be a cyclic-by-$p$ subgroup.  Let $\xi$ be a closed point of $X := \PP^1_k$, and let $\hat Y \to \hat X$ be a connected $I$-Galois cover of $\hat X := \Spec \hat\O_{X,\xi}$.  Then there is a connected $G$-Galois cover $Y \to X: = \PP^1_k$ whose pullback over $\hat X$ is isomorphic to $\Ind_I^G \hat Y$ as a $G$-Galois cover.}

\medskip

\noindent{\it Proof of Lemma~2.5.} By [Ka, Thm.~1.4.1], there is an $I$-Galois cover $f:X_1 \to X$ whose pullback to $\hat X$ is isomorphic to $\hat X_1 \to \hat X$.  In particular, $X_1 \to X$ is totally ramified over $\xi$.  Consider the conjugation action of $I$ on $G$, and form the corresponding semi-direct product $\Gamma = G.I$.  By [Po, Thm.~A], there is a connected $\Gamma$-Galois cover $Z \to X$ that dominates $f:X_1 \to X$, such that $f(B)$ is disjoint from the branch locus of $X_1 \to X$, where $B \subset X_1$ is the branch locus of $Z \to X_1$.  In particular, the inertia group of $Z$ at some point $\zeta$ over $\xi \in X$ is $1.I \subset \Gamma$, and the complete local ring there is isomorphic to that of $X_1$ at the unique point $\xi_1 \in X_1$ over $\xi$.  

Now there is a surjective homomorphism $\Gamma \to G$ given on the first factor of $\Gamma$ by the identity on $G$, and given on the second factor by the inclusion of $I$ into $G$.  The kernel is the normal subgroup $N := \{(i^{-1},i)\,|\,i\in I\} \subset G.I = \Gamma$, which meets $1.I$ trivially.  Let $\phi:Z \to Y := Z/N$ be the corresponding quotient map.  Then $h:Y \to X$ is a connected $G$-Galois cover, whose inertia group at $\eta := h(\zeta)$ is $I \subset G$ (viz.\ the image of $1.I$ under $\Gamma \to \Gamma/N = G$), and whose complete local ring at $\eta$ is isomorphic to that of $Z$ at $\zeta$, or equivalently to that of $X_1$ at $\xi_1$, as an $I$-Galois extension of $\hat\O_{X,\xi}$.  So $Y$ is as desired.  \qed

\medskip

\noindent{\it Remark.} A related result appears as [GS, Theorem 3.4].

\medskip

\noindent{\it Proof of Theorem 2.4.}
 The implication (iii) $\Rightarrow$ (i) is immediate from Theorem~2.2, since each inertia group is a cyclic-by-$p$ subgroup of $G$.  The implication (i) $\Rightarrow$ (ii) is trivial.  So it remains to prove (ii) $\Rightarrow$ (iii).  So let $I = P.C \subset G$ be a cyclic-by-$p$ subgroup of $G$, and let $\hat Y \to \hat X$ be any $I$-Galois cover of $\hat X := \Spec k[[x]]$.  We may identify $\hat X$ with the spectrum of the complete local ring of the affine $k$-line at a point $\xi$.  Applying the lemma, we obtain a connected $G$-Galois cover $Y \to X:=\PP^1_k$ whose pullback to $\hat X$ is $\Ind_I^G \hat Y$.  By (ii), the $G$-Galois cover $Y \to X$ lifts to a $G$-Galois cover ${\cal Y}  \to {\cal X} := \PP^1_R$ for some finite extension $R$ of $W(k)$.  Pulling back to 
the spectrum of $\hat \O_{X,\xi} \approx R[[x]]$, and restricting to 
the identity component of the cover (i.e.\ the component whose closed fibre corresponds to the identity coset of $I$ in $G$), 
we obtain a lifting of $\hat X$ to an $I$-Galois cover $\hat{\cal Y} \to \hat{\cal X} := \Spec R[[x]]$.  This shows that $I$ is a local Oort group, proving (iii).  \qed

\medskip

\noindent{\bf Corollary 2.6.} 
{\sl If a cyclic-by-$p$ group $G$ is an Oort group for $k$, then $G$ is a local Oort group for $k$.}

\medskip

\noindent{\it Proof.} Since $G$ is an Oort group, Theorem~2.4 implies that every cyclic-by-$p$ subgroup of $G$ is a local Oort group of $k$.
In particular, $G$ is a local Oort group of $k$.  \qed

\medskip

\noindent{\bf Corollary 2.7.} 
{\sl If $G$ is an Oort group for $k$, and if $H$ is a subquotient of $G$, then $H$ is an Oort group for $k$.}

\medskip

\noindent{\it Proof.} It suffices to show that every subgroup, and every quotient group, of an Oort group for $k$ is also an Oort group for $k$.  

If $H$ is a subgroup of an Oort group $G$, then every cyclic-by-$p$ subgroup of $G$ is a local Oort group by (i) $\Rightarrow$ (iii) of Theorem~2.4.  In particular, this is the case for every cyclic-by-$p$ subgroup of $H$.  So (iii) $\Rightarrow$ (i) of Theorem~2.4 implies that $H$ is an Oort group. 

If instead $H = G/N$ is a quotient group of $G$, then consider any connected $H$-Galois cover $Y \to \PP^1_k$.  According to the Geometric Shafarevich Conjecture (cf.\ [Ha1], [Po]), the absolute Galois group of the function field $k(x)$ of $\PP^1_k$ is free of infinite rank; so 
there is a connected $G$-Galois (branched) cover $Z \to \PP^1_k$ that dominates $Y \to \PP^1_k$.  Since $G$ is an Oort group, the $G$-Galois cover $Z \to \PP^1_k$ lifts to characteristic $0$, say to ${\cal Z} \to \PP^1_R$.  Let ${\cal Y} = {\cal Z}/N$.   By Proposition~2.1(b),  ${\cal Y} \to \PP^1_R$ is an $H$-Galois cover that lifts $Y \to \PP^1_k$.  This shows that $H$ is an Oort group.  \qed

\medskip  

\noindent{\bf Corollary~2.8.}  
{\sl Let $G$ be a finite group.  Then $G$ is an Oort group if and only if every cyclic-by-$p$ subgroup $I \subset G$ is an Oort group.}

\medskip

\noindent{\it Proof.} The forward implication is immediate from Corollary~2.7.  For the reverse implication, suppose that 
every cyclic-by-$p$ subgroup $I \subset G$ is an Oort group.  Then each such $I$ is a local Oort group, by Corollary~2.6.  So the implication (iii) $\Rightarrow$ (i) of Theorem~2.4 concludes the proof.   \qed

\medskip

 Note that this shows that the latter condition in Corollary~2.8 is equivalent to the three conditions appearing in Theorem~2.4 (i.e.\ we may omit the word ``local'' in (iii) of Theorem~2.4).

Corollary~2.8 reduces the study of Oort groups to the study of cyclic-by-$p$ Oort groups. 

\medskip

\noindent{\bf Proposition 2.9.}  
{\sl Let $n\ge 1$.  If the cyclic group of order $p^n$ is an Oort group for $k$, then so is the cyclic group of order $p^nr$ for every $r$ not divisible by $p$.

Hence the Oort conjecture holds provided that it holds for  cyclic $p$-groups.}

\medskip

\noindent{\it Proof.}  In order to show that $C_{p^nr}$ is an Oort group, it suffices by Theorem~2.4 to show that every subgroup is a local Oort group; each of those is of the form $C_{p^ms}$ for $m\le n$ and $s|r$.  
Suppose $\hat Z \to \hat X := \Spec k[[x]]$ is a connected $C_{p^ms}$-Galois cover.
Let $\hat Y \to \hat X$ be the associated quotient $C_{p^m}$-Galois cover.  Since we assume $C_{p^n}$ is an Oort group, $C_{p^m}$ is a local Oort group by
Corollaries 2.7 and 2.6.  So there is a $C_{p^m}$-Galois cover $\hat{\cal Y} \to \hat{\cal X} := \Spec R[[x]]$ which lifts $\hat Y \to \hat X$ for some discrete valuation ring $R$ which is a finite extension of $W(k)$.
Let $\eta$ be the closed point of $\hat{\cal Y}$ (and of $\hat Y$).  Since $\hat Y$ is smooth over $k$, its lift $\hat{\cal Y}$ is smooth over $R$ and hence regular (at $\eta$).  We want to dominate this by a $C_{p^ms}$-Galois cover that lifts $\hat Z \to \hat X$.

After enlarging $R$, we may assume that each of the codimension $1$ branch points and ramification points of $\hat{\cal Y} \to \hat{\cal X}$ are defined over $R$.  Suppose $\hat{\cal Y} \to \hat{\cal X}$ is not totally ramified over an $R$-point of $\hat{\cal X}$.  Since the subgroups of $C_{p^n}$
are totally ordered, there would then be a proper subgroup $H$ of $C_{p^n}$ such that $\hat{\cal Y}/H \to \hat{\cal X}$
is unramified in codimension $1$.  By purity of the branch locus, $\hat{\cal Y}/H \to \hat{\cal X}$
would then be a non-trivial connected \'etale cover, and hence so would be its special fibre ${\hat Y}/H \to \hat X$.  This is 
impossible by Hensel's Lemma because the residue field $k$ is algebraically closed.  Therefore there
is an $R$-point $P \subset \hat{\cal X}$ that totally ramifies in $\hat{\cal Y}$.   Let $Q \subset \hat{\cal Y}$ be the unique $R$-point over $P$, and let $y \in \hat\O_{\hat{\cal Y},\eta}$ be an element defining the codimension $1$ subscheme $Q$ (which exists since $\hat{\cal Y}$ is regular).  Since 
$\hat\O_{\hat{\cal Y},\eta}$ is complete, it is thus isomorphic to $R[[y]]$.

By Kummer theory (and since $k$ is algebraically closed), the cover $\hat Z \to \hat Y$ is given by $\bar {z}^s=\bar {y}$; here we write $\bar {y} \in \hat\O_{\hat Y,\eta}$ for the residue class of $y \in \hat\O_{\hat{\cal Y},\eta}$ modulo the ideal generated by a uniformizing parameter in $R$.  Let $\hat{\cal Z} \to \hat{\cal Y}$ be the normal $C_s$-Galois cover given by $z^s=y$. Since $s$ is prime to $p$, this is the unique $C_s$-Galois cover of $\hat{\cal Y}$ which lifts $\hat Z \to \hat Y$ and is ramified precisely along $Q$ (again by Kummer theory).  So the composition $\hat{\cal Z} \to \hat{\cal X}$ is Galois, with group $C_{p^ms}$, and it lifts $\hat Z \to \hat X$.  \qed

\medskip

\noindent{\it Remark.} Another approach to Proposition~2.9 would be to 
use that a $C_{p^nr}$-Galois cover is the normalized fibre product of a $C_{p^n}$-Galois cover and a $C_r$-Galois cover.  Namely, if $C_{p^n}$ is an Oort group, then one can lift the unique $C_{p^n}$-Galois quotient cover of a $C_{p^nr}$-Galois cover to a mixed characteristic complete discrete valuation ring $R$; and one can also lift the unique $C_r$-Galois quotient cover using a Kummer extension.  One would then show that if the branch locus of the lift of the $C_r$-Galois cover is chosen suitably (viz.\ as in the above proof), then the normalized fibre product of the two lifts is a smooth cover of $R$-curves, and hence provides the desired lift.  In the cases $n=1,2$, this strategy was carried out explicitly in [GM, II, \S6] by examining equations and relative differents.

\medskip

In the case of {\it local} Oort groups, we have a weaker analog of Corollary~2.7.  First we prove a lemma:

\medskip

\noindent{\bf Lemma~2.10}  
{\sl Let $G=P.C$ be a cyclic-by-$p$ group, with quotient $G'=P'.C'$, where $P,P'$ are $p$-groups and $C,C'$ are cyclic prime-to-$p$ groups.  Then every connected local $G'$-Galois cover $Z' \to X = \Spec k[[x]]$ is dominated by a connected local $G$-Galois cover $Z \to X$, compatibly with the quotient map $G \surj G'$.}

\medskip

\noindent{\it Proof.} Consider the semi-direct product $G''=P'.C$, with $C$ acting on $P'$ through $C'$.  As a first step, we show that $Z' \to X$ is dominated by a $G''$-Galois cover.  Namely, let $Y' \to X$ be the intermediate $C'$-Galois subcover of $Z' \to X$.  By Kummer theory, there is a cyclic $C$ extension $k(W)$ of the function field $k(X) = k((x))$ which 
contains $k(Y')$.  Let $Z''$ be the normalization of $Z'$ in the compositum of $k(W)$ and $k(Z')$ in an algebraic
closure of $k(X)$.  Since $C$ acts on $P' \subset G''$ through $C'$, we have that $Z'' \to X$ is a connected $G''$-Galois cover dominating $Z \to X$.  

To complete the proof, we will dominate the $G''$-Galois cover $Z'' \to X$ by a connected $G$-Galois cover.  Namely, by [Ka, Thm.~1.4.1], $Z'' \to X$ extends to a Katz-Gabber cover, i.e.\ a $G''$-Galois cover $\tilde Y \to \PP^1_k$ whose restriction to $X = \Spec \hat\O_{\PP^1_k,\infty}$ is $Z'' \to X$; whose restriction to $\Spec \hat\O_{\PP^1_k,0}$ is a disjoint union of connected $C$-Galois covers; and which is unramified elsewhere.  Since the kernel of $G \to G''$ is a $p$-group, it follows by [Ha2, Thm.~5.14] (applied to the affine line) that $\tilde Y \to \PP^1_k$ is dominated by a connected $G$-Galois cover $\tilde Z \to \PP^1_k$ such that $\tilde Z \to \tilde Y$ is tamely ramified except possibly over $\infty$ and is \'etale away from $0,\infty$.  

Let $I \subset G$ be an inertia group of $\tilde Z \to \PP^1_k$ over $\infty$.  Since $I$ has $G''$ as a quotient, $I = P''.C$ for some $P''\subset P$.  If $P''$ is a proper subgroup of $P$, then it is contained in a proper normal subgroup $N \subset P$ (since $P$ is a $p$-group); and then $\tilde Z/N$ is an unramified Galois cover of $\tilde Z/P$.  But $\tilde Z/P$ is a $C$-Galois cover of $\PP^1_k$ ramified just at $0,\infty$; hence its genus is $0$ and it has no unramified covers.  This is a contradiction.  So actually $P''=P$, $I= G$, and $\tilde Z \to \tilde X$ is is totally ramified over $\infty$.  (Thus $\tilde Z \to \tilde X$ is a $G$-Galois Katz-Gabber cover.)  Let $\zeta \in \tilde Z$ be the unique point over $\infty$.  Taking $Z= \Spec \hat\O_{\tilde Z, \zeta}$, we have that $Z \to X$ is a connected $G$-Galois cover that dominates $Z'' \to X$ (and hence also $Z' \to X$). \qed

\medskip

\noindent{\bf Proposition~2.11.} 
{\sl If $G$ is a local Oort group for $k$, then every quotient of $G$ is a local Oort group for $k$.}

\medskip

\noindent{\it Proof.} Say $G=P.C$ is an Oort group, with quotient $G'=P'.C'$.  By Lemma~2.10, any connected local $G'$-Galois cover $Z' \to X = \Spec k[[x]]$ is dominated by a connected local $G$-Galois cover $Z \to X$.
Since $G$ is a local Oort group for $k$, the $G$-Galois cover $Z \to X$ lifts to characteristic 0; and taking the corresponding quotient, we obtain a lifting of the given $G'$-Galois cover. \qed

\medskip

We conclude this section with some examples.

\medskip

\noindent {\bf Examples 2.12.}  
As above, $k$ is an algebraically closed field of characteristic $p>0$, and we consider Oort groups and local Oort groups for $k$.

\smallskip

a) Groups of order prime to $p$ are Oort groups for $k$, because all tamely ramified covers lift to characteristic $0$ [Gr, Exp.~XIII, \S2].  {\it Cyclic} prime-to-$p$ groups are also {\it local} Oort groups (e.g.\ by Corollary~2.6, or by [Gr, Exp.~XIII, \S2] applied locally).

\smallskip

b)  By [OSS], the cyclic group $C_p$ is an Oort group, as is $C_{pr}$ with $(p,r)=1$.  By [GM], $C_{p^2}$ is an Oort group, as is $C_{p^2r}$ with $(p,r)=1$.  Since these groups are cyclic-by-$p$ groups, they are also local Oort groups, by Corollary~2.6.  It is unknown whether $C_{p^n}$ is an Oort group for any $n \ge 3$.  

\smallskip

c) It was shown in [GM, I, Example~5.3] that $C_p \times C_p$ is not a local Oort group if $p>2$.  Hence it is also not an Oort group, by Corollary~2.6 above.  By Corollary~2.7 and Proposition~2.11, it then follows that the elementary abelian group $C_p^n$ is neither an Oort group nor a local Oort group for $n>1$ if $p>2$.  Here is a simpler argument, which avoids the machinery of [GM]:  $C_p^n$ acts on the affine line by translation by $\FF_{p^n}$, and hence it acts on the projective line with one fixed point ($\infty$). Taking the quotient by this group, we get a genus $0$ Galois cover of the line in characteristic $p$, with precisely one branch point, where it is totally ramified.  By Corollary 2.3(c), this cover cannot be lifted since $C_p^n$ is not isomorphic to a subgroup of $\PGL_2 = \Aut(\PP^1)$ in characteristic $0$ [Su, Thm.~6.17]; so $C_p^n$ is not an Oort group.  Applying Theorem~2.2 to the above $C_p^n$-Galois cover shows that $C_p^n$ is also not a local Oort group.  (But for every $n$ there exists a local $C_p^n$-cover that lifts [Ma].)

\smallskip

d) For every odd prime $p$, the dihedral group $D_{2p}$ of order $2p$ is a local Oort group [BW, Theorem~1.2].  By Examples (a) and (b) above, every subgroup of $G$ is a local Oort group.  So by Theorem~2.4, $D_{2p}$ is an Oort group.

\smallskip

e) The Klein group $C_2^2$ is an Oort group if $p=2$ (thesis of G.~Pagot [Pa]), and hence a local Oort group.  But $C_2^n$ is {\it not} an Oort group for $n>2$ if $p=2$.  This follows as in Example~(c), since $C_2^n$ acts on the projective $k$-line with one fixed point, but it is not a subgroup of ${\rm PGL}_2 = \Aut(\PP^1)$ in characteristic $0$ [Su, Thm.~6.17].  (On the other hand, $C_2^2$ is a subgroup of ${\rm PGL}_2$ in characteristic $0$.)

\smallskip

f) The quaternion group $Q_8$ of order 8 is not a local Oort group if $p=2$, nor is $\SL(2,3)$.  Namely, the group $\SL(2,3) = Q_8.C_3$ and its subgroup $Q_8$ act faithfully on a supersingular elliptic curve $E$ over $k$, each corresponding to a Katz-Gabber cover of $\PP^1_k$, with the origin of $E$ as the totally ramified point. (In [Si, Appendix~A], see the proof of Prop.~1.2 and Exercise~A.1.)  But $Q_8$ and $\SL(2,3)$ do not act faithfully on any elliptic curve in characteristic $0$; so the assertion follows from Corollary~2.3(b).  By Corollary~2.6, these two groups are also not Oort groups for $k$.

\smallskip

g)  I.~Bouw has announced that the alternating group $A_4 = C_2^2.C_3$ is a local Oort group if $p=2$ (unpublished; see [BW, \S1.3]).  
That implies that every subgroup of $A_4$ is a local Oort group (using Examples (a) and (d) above), and hence that $A_4$ is an Oort group in characteristic $2$, by Theorem~2.4.

\bigskip

\noindent{\bf \S3. Oort groups in odd characteristic.}

\medskip

The main result of this section is that in odd characteristic $p$, every local Oort group, and hence every prime-to-$p$ Oort group, is either a cyclic group $C_n$ or else is a dihedral group of order $2p^n$ for some $n$.  This also has consequences for the structure of arbitrary Oort groups.  We begin with a group-theoretic reduction result:

\medskip

\noindent{\bf Proposition 3.1.} 
{\sl Let $p$ be an odd prime and let $G$ be a finite group 
with a normal Sylow $p$-subgroup $S$ such that $G/S=C$ is cyclic (of order prime to $p$).   
Assume that $G$ has no quotient of the following types:

(1) $C_p \times C_p$;

(2) $P.C_m$, where $P$ is an elementary abelian $p$ group, $p /\!\! | m \ge 3$, and 
$C_m$ acts faithfully and irreducibly on $P$;

(3) $C_p^2.C_2$ where  $C_2$ acts on $P: = C_p^2$
by inversion;

(4) $D_{2p} \times C_\ell$ for some prime number $\ell  > 2$ 
 (including the possibility that $\ell = p$);
 
(5) $C_p.C_4$ where a generator of $C_4$ acts on $P := C_p$ by inversion.

\noindent Then either $G$ is cyclic or it is dihedral of order $2p^n$ for some $n$.}

\medskip

\noindent{\it Proof.}   We proceed inductively, and we assume that the proposition holds for every group of order less than $\# G$.  
We may assume that $p$ divides the order of $G$ (for otherwise
$G \approx C$ is cyclic).  Since $S$ and $G/S$ have relatively prime orders, $G$ contains a subgroup
isomorphic to $G/S$, which we again denote by $C$.  Set $K=C_C(S)$.  Then every subgroup of $K$ is normal in the cyclic group $C$
and is normalized by $S$, and hence is in normal in $G$.  In particular, $K$ is normal in $G$.

Suppose $G$ has the property that it has no quotient of the form (1)-(5) and that $H$ is a quotient of $G$.
Then $H$ has a normal Sylow $p$-subgroup and the quotient of $H$ by this subgroup is cyclic of order prime to $p$.
Furthermore, $H$ can have no quotient of the form (1) - (5).     So by the inductive hypothesis, every proper quotient $H$ of $G$ is cyclic or is dihedral of order $2p^m$ for some $m$.    
Hence if $N$ is any non-trivial  normal subgroup
of $G$ contained in $S$, then $G/N$ is either cyclic or else dihedral of order $2p^n$.
In particular, this implies that $S/N$ is cyclic.  

The Frattini subgroup $\Phi(S)$ of $S$ is normal in $G$
since $S$ is normal.  Suppose $\Phi(S)$ is non-trivial.  Then $S/\Phi(S)$ is cyclic, so $S$ is cyclic by the Burnside Basis Theorem.  If $G = S$ then $G$ is cyclic.  If $G \ne S$  then $G/\Phi(S)$ is a proper
quotient of $G$ that is not a $p$-group but which has order divisible by $p$ since $\Phi(S) \ne S$.  Hence
$G/\Phi(S)$ is dihedral of order $2p^m$ for some $m>0$ and $\#(G/S) = 2$.  In this case an involution in $G$ either centralizes $S$ (and so $G$ is cyclic) or acts as inversion on $S$ (and so $G$ is dihedral).
This completes the proof if $\Phi(S)$ is non-trivial.

We now suppose that $\Phi(S)$ is trivial or equivalently
that $S$ is elementary abelian.  If $S$ is central in $G$, then $G \approx S \times C$, and
$G$ surjects onto $S$.  Since $G$ does not surject onto $C_p \times C_p$, neither does $S$.  So the elementary abelian $p$-group $S$ is isomorphic to $C_p$, and hence $G  \approx S \times C$ is cyclic.  So from now on we may assume that $S$ is not central, i.e.\
$C$ does not commute with $S$.  Thus 
a generator $x$ for $C$ induces an automorphism
of order $m > 1$ on $S$, by conjugation.   

Consider the case that $m > 2$.  Then the action of $C$ on $S$ cannot be both faithful and irreducible, since then $G$ would be a group as in (2), a contradiction.  On the other hand, if $C$ does not act faithfully on $S$, then the normal subgroup $K=C_C(S)$ is non-trivial; and hence the quotient $G/K$ is either cyclic or dihedral, which contradicts the assumption that $m>2$.  Finally, suppose $C$ does not act irreducibly on $S$.  Since $S$ is an elementary abelian $p$-group and $C$ is cyclic of order prime to $p$, $S$ is the product $\prod_{i = 1}^t S_i$ of some number $t > 1$
of subgroups $S_i$ on which $C$ acts irreducibly by conjugation.  For each $1 \le j \le t$, $T_j = \prod_{i \ne j} S_i$
is a non-trivial normal subgroup of $G$, so $G/T_j$ is either cyclic or dihedral.  This means $C$ acts
trivially or by inversion on $S/T_j \approx S_j$ for all $j$, which contradicts the assumption that $m > 2$. 

Therefore $m=2$.   Suppose that the elementary abelian $p$-group $S$ is not cyclic.  Then there exists a $C$-invariant subgroup $T$ of $S$ having index $p^2$.  Since $G/T$ contains a subgroup $S/T$ that is isomorphic to $C_p \times C_p$, it is neither cyclic nor dihedral.  It follows that the normal subgroup $T$ is trivial and so $S \approx C_p^2$.  With $K=C_C(S)$ as above, since $m=2$ we have that $G/K$ is either of the form (3) or (4) in the statement of the result, with $\ell=p$ in the case of (4).  This is a contradiction.  

So we are reduced to the case that $S$ is cyclic of order $p$, and $m=2$.  If $K=C_C(S)$ is non-trivial, let $K'$ be a maximal proper subgroup of $K$.  Then $K'$ is normal in $G$, and $G/K'$ is of the form $C_p.C_{2\ell}$ for some prime $\ell$, where the generator of $C_{2\ell}$ acts by inversion.  Depending on whether $\ell$ is odd or is equal to $2$, $G/K'$ is then of the form (4) (with $\ell \ne p$) or (5).  This is a contradiction.  So in fact $K$ is trivial, hence 
$G$ is dihedral of order $2p$. \qed

\medskip

In order to apply Proposition~3.1, we show in the next result that certain groups are not local Oort groups.  In the proof, we use that for any polynomial $f(u)$ of degree $m$ prime to $p$, the genus of the characteristic $p$ curve $w^p-w=f(u)$ is $(p-1)(m-1)/2$.  This formula follows from the tame Riemann-Hurwitz formula, viewing the curve as a cover of the $w$-line.

\medskip

\noindent{\bf Proposition~3.2.} 
{\sl The groups listed in items (1)-(5) of Proposition~3.1 are not local Oort groups for an algebraically closed field $k$ of odd characteristic $p$.}

\medskip

\noindent{\it Proof.}  The case of type (1) of Proposition~3.1 was shown in [GM, I, Example 5.3]; see also Example 2.12(c) above.  So it remains to consider types (2)-(5).

\smallskip

In types (2) and (3), $G$ is isomorphic to a subgroup of $\PGL(2,k)$ consisting of upper triangular matrices, by [Su, Thm.~6.17].  So we obtain an action of $G$ on $Y := \PP^1_k$ such that the $G$-Galois cover $Y \to X = Y/G$ is totally ramified at infinity and only tamely ramified elsewhere.  Here $X$ necessarily has genus $0$; so we obtain a genus $0$ Katz-Gabber $G$-Galois cover of $X := \PP^1_k$.  But $G$ cannot be embedded into ${\rm PGL}_2(K)$ for any field $K$ of characteristic $0$ [Su, Thm.~6.17].  So by Corollary 2.3(c), the local cover $\hat Y \to \hat X$ obtained by completing $Y \to X$ at infinity cannot lift to characteristic $0$.

\smallskip

In type (4), first consider the situation of $\ell = p$.  
Let $X$ be the projective $x$-line over $k$ and let 
$Y \to X$ be the $G$-Galois Katz-Gabber cover given by $t^2=x$, $u^p-u=t$, $v^p-v=x$.  This cover is totally ramified over $x=\infty$ and tamely ramified of index $2$ over $x=0$ (and unramified elsewhere).   Rewriting the equations by eliminating $x$ and $t$, the curve $Y$ is given by by the equation $v^p-v=(u^p-u)^2$; or equivalently by $w^p-w=-2u^{p+1}+2u^2$ (setting $w=v-u^2$).  Applying the genus formula given just before the statement of the proposition, we find that the genus of $Y$ is $p(p-1)/2$. Let $T \to X$ be the quotient cover of $Y \to X$ given by $t^2=x$ and let $H = \Gal(Y/T)$.  So $T \to X$ is a degree $2$ tame cover of genus $0$, branched at two points.

Now suppose that there is a curve $Y^\circ$ of genus $p(p-1)/2$ in characteristic $0$ and a faithful action of $G$ on $Y^\circ$ whose quotient $X^\circ := Y^\circ/G$ has genus $0$.  Let $T^\circ = Y^\circ/H$. So $T^\circ \to X^\circ$ is a degree $2$ cover of genus $0$ (since the genus of $T$ is $0$), and hence $T^\circ \to X^\circ$ is branched at two points.  Also, $Y^\circ \to T^\circ$ is a $C_p^2$-Galois cover, say with $n$ branch points; here $n>2$ since the cover $Y^\circ \to T^\circ$ is not cyclic.  So over each of these $n$ branch points, $Y^\circ \to T^\circ$ has $p$ ramification points, each with ramification index $p$.  By the characteristic $0$ Riemann-Hurwitz formula, we have that $p(p-1)-2=-2p^2 + np(p-1)$.
Rearranging and dividing by $p-1$ gives $2(p+1)=(n-1)p$, which is impossible since the odd prime $p$ does not divide the left hand side.  So in this case the result follows from Corollary~2.3(b).

\smallskip

It now remains to consider the case  in which $G$ is of type (4) with $\ell \ne p$ or of type (5). Then $G$ is  a semi-direct product $C_p . C_{2\ell}$ with a generator of $C_{2\ell}$ acting by inversion on $C_p$. Let $T$, $X$, $Y$ and $Z$ be copies of the
projective line $\PP^1_k$ with affine coordinates $t$, $x$, $y$ and $z$, respectively.  Define cyclic covers $X \to T$,
$Y \to X$ and $Z \to X$ with groups $C_2$, $C_\ell$ and $C_p$, respectively, by $ t= x^2$, $x = y^\ell$ and $x = z^p - z$.
Then $Y \to T$ is defined by $t = y^{2\ell}$ and is a Katz-Gabber $C_{2\ell}$-Galois 
cover, while $Z \to T$ is a Katz-Gabber $D_{2p}$-Galois cover.  We find that if $W$ is the normalization of $Z \times_X Y$, then
$W \to T$ is a Katz-Gabber $G$-Galois cover. Since $W \to Y$ is defined by $z^p - z = y^\ell$, the formula in the paragraph
just prior to the statement of Proposition 3.2 shows that $W$ has genus $g_W = (p-1)(\ell-1)/2$.

Suppose now that the $G$-Galois cover $W \to T$ lifts to characteristic $0$.    By taking the base change of such a lift to an algebraically closed field $L$, we obtain a $G$-Galois cover $W^\circ \to T^\circ$  of $L$-curves with the following properties.  By Proposition 2.1, $g_W = g_{W^\circ}$ and the curves $Z^\circ = W^\circ/C_{\ell}$, 
$Y^\circ = W^\circ/C_{p}$ and $T^\circ = W^\circ/G$ have genus $0$ since this is true of the corresponding quotients of $W$.  
Since $L$ is algebraically closed, $Y^\circ$
is isomorphic to $\PP^1_L$.  Because ${\rm char}(L)=0$, each non-trivial element of ${\rm Aut}(\PP^1_L) = {\rm PGL}_2(L)$
of finite order is conjugate to the class of a diagonal matrix, and thus fixes exactly two points of $\PP^1_L$.   Hence the branch locus of the $C_{2\ell}$-Galois cover $\pi_{Y^\circ}:Y^\circ \to T^\circ$ consists of
two totally ramified points $\{Q_1,Q_2\} \subset T^\circ$.   The inertia group in $G$ of each point of $W^\circ$ over
$Q_i$ is cyclic, since ${\rm char}(L) = 0$, and of order divisible by $2\ell$.  So these inertia groups have order
$2\ell$.  There are now $2(\#G)/(2\ell) = 2p$ points over $\{Q_1,Q_2\}$ in $W^\circ$, which all ramify in the  tame $C_\ell$-Galois cover $\pi:W^\circ \to Z^\circ$ as $C_\ell$ is normal in $G$.  The Riemann-Hurwitz formula for $\pi$ now gives 
$$ g_{W^\circ} \ge 1 + \ell(g_{Z^\circ} - 1) + p(\ell -1) = (-1 + p)(\ell - 1) > (p-1)(\ell -1)/2 = g_W$$
since $g_{Z^\circ} = 0$ and $g_W > 0$.  This contradicts $g_{W^\circ} = g_W$, which completes the proof.\qed

\medskip

As a consequence, we obtain

\medskip

\noindent{\bf Theorem 3.3.} 
{\sl  Suppose that $p = {\rm char}\,k > 2$ and $G$ is a local Oort group for $k$.  Then $G$ is either cyclic or is isomorphic to a dihedral group of order $2p^n$ for some $n$.}

\medskip

\noindent{\it Proof.} If $G$ is a local Oort group for $k$, then so is every quotient of $G$, by Proposition~2.11.  So by Proposition~3.2, the groups listed as items (1)-(5) in Proposition~3.1 cannot be quotients of $G$.  Thus by Proposition~3.1, $G$ is of the asserted form.  \qed

\medskip

By Corollary 2.6, this theorem implies the forward direction of the Strong Oort Conjecture in odd characteristic $p$:

\medskip

\noindent{\bf Corollary 3.4.} 
{\sl Suppose that $p = {\rm char}\,k > 2$ and $G$ is a cyclic-by-$p$ group.  If $G$ is an Oort group for $k$, then $G$ is isomorphic to some $C_n$ or $D_{2p^n}$.}

\medskip

By Corollary~2.8 this in turn implies

\medskip

\noindent{\bf Corollary 3.5.} 
{\sl If $G$ is an Oort group in odd characteristic $p$, then every cyclic-by-$p$ subgroup of $G$ is isomorphic to some $C_n$ or $D_{2p^n}$.}

\medskip

The consequences of Corollary~3.5 will be explored further in [CGH2].  For now we note these corollaries of the above results:

\medskip

\noindent{\bf Corollary 3.6} 
{\sl Let $G$ be an Oort group for $k$, where $k$ has characteristic $p>2$.  Then the Sylow $p$-subgroups of $G$ are cyclic.}

\medskip

\noindent{\it Proof.}  Let $P$ by a Sylow $p$-subgroup of $G$.  By Corollary~2.7, the subgroup $P \subset G$ is an Oort group for $k$.  
Since $p \ne 2$, a dihedral group $D_{2p^n}$ is not a $p$-group, and so is not isomorphic to $P$.  So Corollary~3.4 implies that $P$ is cyclic.  \qed

\medskip
 \noindent{\bf Corollary 3.7.} 
    {\sl Let $G$ be an Oort group for $k$, where ${\rm char}\,k = p > 2$.
    Let $P \subset G$ be a $p$-subgroup of $G$, and suppose that some $g \in
    G$ normalizes $P$ but does not centralize $P$.  Then $g$ has order $2$,
    and $g$ acts by inversion on $P$ and on the abelian subgroup $Z :=
    C_G(P)$, where $Z$ is also equal to $C_G(S)$ for any  Sylow $p$-subgroup
    $S$ containing $P$.}

    \medskip

    \noindent{\it Proof.} Let $C$ be the subgroup generated by $g$.  Since
    $g$ normalizes $P$, the subgroup $I$ generated by $P$ and $g$ is
    a cyclic-by-$p$ subgroup of $G$.  By Corollary~2.7, $I$ is an Oort group
    for $k$.  By Corollary~3.4, $I$ either is cyclic or is dihedral of order
    $2p^n$.  The former case is impossible because $g$ is assumed not to
    centralize $P$.  The latter case implies that $g$ has order $2$ and that
    the conjugation action of $g$ on $P$ takes each element to its inverse.

    If $z \in Z$, then $gz$ normalizes but does not centralize $P$.  So by
    the previous paragraph, $gz$ is an involution.  Since $g$ is also an
    involution, $gzg^{-1}=z^{-1}$; i.e.\ $g$ acts by inversion on $Z$.
    Since inversion is an automorphism of $Z$, $Z$ is abelian.  Now $P
    \subset S$, so $C_G(S) \subset C_G(P) = Z$.  But $S$ is abelian by
    Corollary~3.6, so $S \subset C_G(S) \subset Z$ and hence $C_G(Z) \subset
    C_G(S)$.  Since $Z$ is also abelian, $Z \subset C_G(Z) \subset C_G(S)
    \subset Z$, i.e.\ all these groups are equal.
    \qed
 
\bigskip
\noindent{\bf \S4.  Oort groups in characteristic two.}

\medskip

The classification of Oort groups in characteristic two is more involved than in odd characteristic.  In this section we show that a cyclic-by-$2$ Oort group in characteristic $2$ is either cyclic, or a dihedral $2$-group, or is the alternating group $A_4$.  We also show a corresponding result for local Oort groups.  

We begin by recalling some notation and facts about $2$-groups.
A {\it generalized quaternion group} of order $2^a$, $a \ge 3$, is given by
$Q_{2^a} = \langle x, y | x^{2^{a-1}}=1, yxy^{-1}=x^{-1}, y^2=x^{2^{a-2}} \rangle$.
It follows from [Go,~Chap. 5,  Thm. 4.10(ii)] that these are the only noncyclic $2$-groups that contain a unique involution.  The group $Q_8$ is the usual quaternion group of order $8$.

The {\it semidihedral group} of order $2^a$, $a > 3$, is denoted ${\rm SD}_{2^a}$
and has presentation 
$\langle x, y | x^{2^{a-1}}=1, y^2=1, yxy= x^{-1 + 2^{a-2}} \rangle$.
Note that if $G$ is dihedral, semidihedral or generalized quaternion
then $G/[G,G]$ is elementary abelian of order $4$.  The next lemma
shows that these groups are characterized by this property.

\medskip

\noindent{\bf Lemma 4.1} 
{\sl  Let $G$ be a finite $2$-group
with derived group $D$ and whose abelianization $G/D$ is a Klein four group.
Then $G$ is dihedral, semidihedral or generalized quaternion.
If in addition, $\Aut(G)$ is not a $2$-group, then $G$ is either
a Klein four group or is quaternion of order $8$.}

\medskip

\noindent{\it Proof.}  The first assertion is contained in Theorem 4.5 of Chapter 5 of [Go]. 
For the second assertion, suppose that $\Aut(G)$ is not a $2$-group (and still assume that $G$ is non-abelian).  Then 
$G$ admits an automorphism $\sigma$ of odd order, which necessarily
acts faithfully on $G/D$ and on $G/Z(G)$.  Since $G/D$ is a Klein four group, $\sigma$ has order $3$.  If $\# G \ge 16$, 
then in all cases $G/Z(G)$ is a dihedral $2$-group, with a unique cyclic subgroup of index $2$, which must be invariant under $\sigma$.  This is impossible since $\sigma$ has order $3$.  So actually  $\# G = 8$. The argument just given shows that $G \not\approx D_8$ (since otherwise an automorphism of order $3$ would have to fix the unique cyclic subgroup of index $2$), and so $G \approx Q_8$. \qed

\medskip

Using this lemma, we obtain the following group-theoretic reduction result, which is analogous to Proposition~3.1:

\medskip

\noindent{\bf Proposition 4.2} 
{\sl Let $G$ be a finite group with a normal Sylow $2$-subgroup $S$ such that
$G/S=C$ is cyclic (of odd order).   Assume that $G$ has no
quotient of the following types:

(1)  $P.C$, where $P$ is an elementary abelian $2$-group and $C$ is a cyclic group of odd order at
least $5$ that acts irreducibly on $P$;

(2) $C_2^4.C_3$, where $C_3$ acts
without fixed points on $P := C_2^4$;

(3) $C_4^2.C_3$ where $C_3$ acts
faithfully on $P :=C_4^2$;

(4) $C_2^3.C$, where $C$ has order $1$ or $3$ and acts 
faithfully on $E := C_2^3$ (i.e.\ $G$ is isomorphic to $C_2^3$ or $A_4 \times C_2$); 

(5) $C_2^2 \times C_\ell$ for some odd prime $\ell$; 

(6)  $C_2^2.C_{3\ell}$ where $\ell$ is an odd prime and $C:=C_{3\ell}$ acts nontrivially on $P := C_2^2$;

(7) $C_4 \times C_2$.

\noindent Then $G$ either is a cyclic group, or is isomorphic to $A_4$ or $SL(2,3)$, or is a $2$-group that is dihedral, semidihedral
or generalized quaternion.}

\medskip
 
\noindent{\it Proof.}  As in Proposition~3.1, we proceed inductively by assuming that the proposition holds for every group of order less than $\# G$.  Since $G = S.C$, we may view $C$ as a subgroup of $G$.
Let $K=C_C(S)$ and note that every subgroup of $K$ is normal in $G$ (as in the proof of Proposition~3.1).
By the inductive hypothesis, every nontrivial quotient of $G$ satisfies the conclusion of the theorem.  
We consider various cases for $S/\Phi(S)$, where $\Phi(S)$ is the Frattini subgroup of $S$.

\smallskip

Case 1.  If $S/\Phi(S)$ is cyclic, then so is $S$, by the Burnside Basis Theorem.  Since every automorphism of $S$ has $2$-power order, the odd-order cyclic group $C$ acts trivially on $S$.  Thus $G = S \times C$, which is cyclic.

\smallskip

Case 2.  If $S/\Phi(S)$ has order greater than $4$,
then $\Phi(S)=1$  since otherwise, 
$G/\Phi(S)$ would be a counterexample to the result, contradicting the inductive hypothesis.  Thus in this case, $S$ is elementary abelian of order at least $8$.
Similarly, $K=1$, so $C$ acts faithfully on $S$.   Let $T$ be a non-trivial minimal normal subgroup
of $G$ contained in $S$.  Then the quotient $G/T$ satisfies the hypotheses of the Proposition and has order less that $
\#G$, so by
the inductive hypothesis it has
one of the asserted forms.  Since its Sylow $2$-subgroup is elementary abelian, $G/T$ must be $A_4$, cyclic, or $D_4 = C_2^2$.  If 
$T=S$, which is an elementary abelian $2$-group of order at least $8$, then $C$ acts irreducibly on $S$ since $T$  is minimal.  In this
case $C$ has order at least $5$, so
$G$ is a group as in (1), which is a contradiction.
Alternatively, if $T$ is strictly contained in $S$, then $S=T \times U$ with $U$ normal in $G$ (by complete
reducibility).  Since $\#S \ge 8$, and since $\# T \le \# U$ by minimality of $T$, we have that the elementary abelian $2$-group $U$ has order at least $4$.  So $U.C = G/T$ cannot be cyclic, and hence must be isomorphic to $C_2^2$ or $A_4$.  Thus $U$ is a Klein four group; $\#C=1$ or $3$; and the order of $T$ is $2$ or $4$.  If the normal subgroup $T$ has order $2$, then it is central in $G$; so $G$ is of type (4), a contradiction.  If $T$ has order $4$, then $C$ acts irreducibly on $T$ (by minimality of $T$); so $G$ is of type (2), again a contradiction.

\smallskip

Case 3.  The remaining case is when the elementary abelian $2$-group $S/\Phi(S)$ has order $4$ (i.e.\ is a Klein four group).
We further subdivide this case. 

\smallskip

Case 3(a).   $S$ is nonabelian.   Let $D$ be the derived
subgroup of $S$.   Then $G/D$ is a proper quotient of $G$; and so by the inductive hypothesis, it satisfies the conclusion
of the proposition.  Since $G/D$ has an elementary abelian Sylow $2$-subgroup, $G/D \approx A_4$ or is a Klein four group.
In particular, $S/D$ a Klein four group.  So
by Lemma~4.1,  $S$ is dihedral, generalized
quaternion or semidihedral.   Moreover, $C = G/S$ has order at most
$3$ since $G/D$ is isomorphic to $A_4$ or $C_2^2$.  So $C_C(S/D)=1$ in $G/D$; hence $K=1$ and $C$ acts faithfully on $S$.  If $C$ is trivial, then $G=S$ is a $2$-group of rank $2$ with no quotient isomorphic to $C_4 \times C_2$ (by type (7) of the assertion); hence $G$ is dihedral, semidihedral, or generalized quaternion and the result holds.  The remaining possibility is that $C$ has order $3$ and $S$ admits an automorphism of order $3$,
whence $S$ is quaternion of order $8$ and $G=SL(2,3)$. So again the result holds for $G$.

\smallskip

Case 3(b).   $S$ is abelian, necessarily of order at least $4$.
If $\#S > 4$, then $K=1$ by the inductive hypothesis applied to $G/K$.  So $C$ acts faithfully
on $S$ and hence on $S/\Phi(S)$; thus $\# C  \le 3$.   If $C=1$,
then $G$ surjects onto $C_4 \times C_2$, contradicting (7).
If $C$ has order $3$, then by modding out by the subgroup
generated by $\{s^4\,|\, s \in S\}$, we may assume that $S$
has exponent $4$.  So $S$ is either $C_4 \times C_4$ or
$C_4 \times C_2$.   The first case cannot occur because of
(3) and the second case cannot occur because that group
has no automorphisms of order $3$.  This is a contradiction.

So actually $\# S =4$.  Thus $\Phi(S)=1$ and $S$ is a Klein four group.  If $K \ne 1$, then the inductive hypothesis implies that $K$ has prime order $\ell$, and $G = S \times K$ or $G/K = A_4$ (depending on whether the image of $C$ in $\Aut(S)$ is $1$ or $3$).  But this is impossible, because the former group is ruled out by (5) and the latter group by (6).  So actually $K = 1$, and $C$ acts faithfully on $S$.  Hence $G$ is isomorphic to $C_2^2 = D_4$ or $A_4$ and the result holds for $G$.  
\qed

\medskip

Analogously to Proposition~3.2, we have

\medskip

\noindent{\bf Proposition~4.3} 
{\sl Let $k$ be algebraically closed of characteristic $2$.  Then none of the groups of type (1)-(7) in Proposition~4.2 are local Oort groups for $k$, nor are $Q_8$ and $SL(2,3)$.} 

\medskip

\noindent{\it Proof.} We consider each of these types of groups in turn.
\smallskip

If $G$ is of type (1) or (2), or the first case of type (4), then $G$ embeds into the upper triangular matrices of $\PGL_2(k)$, but it does not embed into $\PGL_2(K)$ for $K$ of characteristic $0$ [Su, Thm.~6.17].  So $G$ is the Galois group of a Katz-Gabber cover of genus $0$ over $k$, but it is not  a local Oort group by Corollary 2.3(c).  

\smallskip

For the next several types of groups, we let $Z$ and $X$ be copies of the projective line over $k$, with affine parameters $z$ and $x$.  Let $\FF_4$ be the field with four elements, and fix an isomorphism of the additive group of
$\FF_4$ with $C_2^2$.
We consider the $C_2^2$-Galois cover $Z \to X$ given by $x=z^4-z$, with $\alpha \in \FF_4 = C_2^2$ acting
on $Z$ by $z \mapsto z+\alpha$.    
Let $t=x^3$, so that $X \to T$ is a $C_3$-Galois cover branched at $t=0, \infty$, where $T$ is the $t$-line.  Then the composition $Z \to X \to T$ is a Galois cover with group $A_4 = C_2^2.C_3$.  Note that $Z \to T$ is a Katz-Gabber cover, and that $Z$ has genus $0$.   The Sylow $2$-subgroup $P_H$ of $H = {\rm Gal}(Z/T)$ is $C_2^2 = \FF_4$.  An
element $\alpha \in \FF_4 = P_H$ sends the uniformizer $z^{-1}$ at the unique point $\infty_Z$ of $Z$ over $t = \infty$
to $(z + \alpha)^{-1} = z^{-1}(1+ \alpha z^{-1})^{-1}$.  We see from this that the second lower ramification group $H_2$ associated to $\infty_Z$ is trivial.  Thus the lower ramification group $(P_H)_v$ is trivial if $v > 1$.  Since $P_H = (P_H)_0 = (P_H)_1$, this implies that the upper ramification group $P_H^u$ is trivial for $u > 1$.
\smallskip

Suppose now that $G$ is  either a group of type (3) or a group of type (4) for which $3|\# G$.  (We already 
treated above the case of groups of type (4) for which $3$ does not divide $\# G$.) Then $G$ is an extension of $A_4$ by
a minimal normal group $N$ isomorphic to either $C_2$ or $C_2^2$.  Identify $H$ with $G/N \approx A_4$.  
By Lemma~2.10 there exists a local $G$-Galois cover
dominating the completion of the above $H$-Galois cover $Z \to T$ at its totally ramified point; and so there is also a corresponding $G$-Galois Katz-Gabber cover $S \to T$ dominating $Z \to T$.  Let
$g$ denote the genus of $S$ and let $P$ be the Sylow $2$-subgroup of $G$.  We consider the upper and lower ramification groups of $P$ and $P_H$ at the totally ramified points of $S \to T$ and $Z \to T$.   
By [Se,  \S IV.3, Prop. 14], $(P^u\cdot N)/N = (P/N)^u = (P_H)^u$ for all $u$, and we have shown this is trivial for $u > 1$.  Thus
$P^u \subset N$ for $u > 1$.  Since $N$ is a minimal normal subgroup of $G$, $P^u$ is either equal to
$N$ or trivial for $u > 1$.  Since $P_0 = P_1$, this implies $P_v$ is either $N$ or trivial for $v > 1$.  By the Hasse-Arf theorem [Se, IV, \S3], the number of $i$ such that $P_i=N$  is divisible by $(P:N)=4$. It follows that sequence of lower ramification groups of $P$ has the form
$P = P_0 = P_1$, $N = P_2 = \cdots = P_{1 + 4a}$ for some $0 \le a \in \ZZ$, and $P_i = \{e\}$ for $i \ge 2 + 4a$.
Thus the wild form of the Riemann-Hurwitz formula
for the $N$-Galois cover $S \to Z = S/N$ gives
$$2g_S - 2 = \#N (2g_Z - 2) + (\#N - 1)(2+4a) \equiv 2 \quad {\rm mod} \quad 4$$
using $g_Z = 0$. Thus $g_S$ is even.   Let $D$ be the set of orders of non-trivial elements of $G$.
Suppose $S^\circ \to T^\circ$ is a $G$-Galois cover of smooth connected curves over an algebraically
closed field of characteristic $0$ and that $T^\circ$ has genus $0$.  The tame Riemann-Hurwitz formula shows 
$$
2(g_{S^\circ}-1) =  \#G \left (-2 + \sum_{d \in D} b_d (d-1)/ d \right ) \equiv 0 \quad {\rm mod}\quad 4,
$$
where $b_d$ is the number of branch points with (cyclic) inertia groups of order $d$ and $\# G/d \equiv 0 \equiv 2\# G$ mod $4$
for $d \in D$.  So $g_{S^\circ}$ is odd and cannot equal $g_S$. Thus, the above Katz-Gabber cover $S \to T$ cannot lift to characteristic $0$.  So by Corollary 2.3(a), this completes the proof that no group of type  (3) or (4) can be a local Oort group.

\smallskip

To treat $G$ as in cases (5) and (6), we first construct a $C_2^2.C_{3\ell}$-Galois cover $V \to T$ of the projective line $T$
over $k$. Let $Z\to T$ be the $A_4 = C_2^2.C_3$-Galois cover constructed previously, with quotient
$C_3$-Galois cover $X \to T$ defined by $t = x^3$ on affine coordinates for the projective lines $X$ and $T$,
respectively.  Let $Y\to T$ be the $C_{3\ell}$-Galois cover of projective lines defined on affine coordinates
by $t = y^{3\ell}$. This has subcover $Y \to X$ defined by $x = y^\ell$.  The normalization $V$
of the fibre product $Z \times_X Y$ now gives a $C_2^2.C_{3\ell}$-Galois cover $V \to T$ which has
a $C_2 \times C_\ell$-Galois subcover $V \to X$.  It will suffice to show that this subcover cannot be lifted
to characteristic $0$.  If there were such a lift, then after making a a base change to an algebraically
closed field $L$ of characteristic $0$ we would have  a $C_2^2\times C_\ell$-Galois cover $V^\circ \to X^\circ$
of smooth connected projective curves over $L$ such that $g_{V^\circ} = g_V$, $g_{X^\circ} = g_X = 0$,
$g_{Z^\circ} = g_Z = 0$ when $Z^\circ = V^\circ/C_\ell$  and $g_{Y^\circ} = g_Y = 0$
when $Y^\circ = V^\circ/C_2^2$.  Since $Z^\circ \to X^\circ$ and $Y^\circ \to X^\circ$
have groups $C_2^2$ and $C_\ell$ of coprime orders and $G = C_2 \times C_\ell$, the 
branch locus $B^\circ$ of the $C_\ell$-Galois cover $V^\circ \to Z^\circ$ is the pullback via $Z^\circ \to X^\circ$
of the branch locus of $Y^\circ \to X^\circ$.  Thus $B^\circ$ is taken to itself by the action of $C_2^2 = {\rm Gal}(Z^\circ/X^\circ)$, so since
inertia groups in characteristic $0$ are cyclic we see that $\# B^\circ$ is even.  However,
the same argument shows that the branch locus $B$ of $V \to Z$ is the pullback via
$Z \to X$ of the branch locus $x \in \{0,\infty\}$ of $Y \to X$.  Since $Z \to X$ was defined
by the affine equation $z^4 - z = x$, we see that $\# B = 5$, so $\# B^\circ \ne \#B$.
However, this contradicts $g_Z = g_{Z^\circ}$, $g_{V} = g_{V^\circ}$ and the tame
Riemann-Hurwitz formulas for the $C_\ell$-Galois covers $V \to Z$ and $V^\circ \to Z^\circ$.  The 
contradiction completes the treatment of cases (5) and (6).

\smallskip

The group $G=C_4 \times C_2$, of type (7), acts on the genus $2$ curve $X: y^2-y=x^5$ in characteristic $2$, with commuting generators $\sigma$, $\tau$, of orders $4$, $2$ respectively, given by $\sigma(x,y)=(x+\zeta,y+\zeta^2x^2+\zeta x + \xi)$, $\tau(x,y)=(x+1,y+x^2+x+\zeta)$, where $\zeta$ is a primitive cube root of unity and $\xi^2-\xi = \zeta^2$.  The quotient morphism $X \to X/G$ is a $G$-Galois cover with a unique ramification point (the point at infinity), which is totally ramified.  By the wild form of the Riemann-Hurwitz formula, $X/G$ has genus $0$; i.e.\ this is a Katz-Gabber cover of the line with group $G$.  But by the tame Riemann-Hurwitz formula, any $G$-Galois cover of the line in characteristic $0$ must have odd genus (using that the number of branch points with ramification index $4$ must be even).  So the Katz-Gabber cover cannot lift to characteristic $0$, and Corollary~2.3(a) implies that $C_4 \times C_2$ is not a local Oort group.    

\smallskip

The last assertion is contained in Example~2.12(f).
\qed

\medskip

\noindent{\it Remark.} For groups $G$ of type (5) in the above result, even more is true: {\it no} local $G$-Galois covers lift to characteristic zero.  This follows from a result of Green and Matignon [Grn, Cor.~3.3], saying that for an abelian cover to lift, the group must be cyclic or a $p$-group.  

\medskip

The above results yield the following analogs of Theorem~3.3 and its corollaries:

\medskip

\noindent{\bf Theorem~4.4.} 
{\sl   Suppose that ${\rm char}\,k = 2$ and $G$ is a local Oort group for $k$.  Then $G$ is either cyclic, or is isomorphic to a dihedral $2$-group, or is isomorphic to $A_4$, or is isomorphic to a semi-dihedral group or generalized quaternion group of order $\ge 16$.}

\medskip

\noindent{\it Proof.} By Proposition~2.11, every quotient of $G$ is also a local Oort group for $k$.  So by Proposition~4.3, $G$ is not isomorphic to $\SL(2,3)$ or $Q_8$, and no quotient of $G$ is isomorphic to a group of type (1)-(7) in the statement of Proposition~4.2.  That latter proposition then implies the theorem. \qed

\medskip

\noindent{\it Remark.}   In [CGH1], we will show that in fact semi-dihedral groups are not local Oort groups in characteristic $2$; the status of generalized quaternion groups  of order $\ge 16$ as local Oort groups remains open.  See also the remark after Theorem~4.5.

\medskip

The following theorem is the forward direction of the Strong Oort Conjecture in characteristic $2$.

\medskip

\noindent{\bf Theorem~4.5.} 
{\sl Suppose that ${\rm char}\,k = 2$ and $G$ is a cyclic-by-$2$ group.  If $G$ is an Oort group for $k$, then $G$ is either cyclic, or is isomorphic to a dihedral $2$-group, or is isomorphic to $A_4$.}

\medskip

\noindent{\it Proof.} If $G$ is a cyclic-by-$2$ Oort group for $k$, then $G$ is a local Oort group of $k$ by Corollary~2.6.  Hence $G$ is one of the possibilities listed in Theorem~4.4. By  Corollary~2.7, every subgroup of $G$ is also an Oort group.  But the quaternion group $Q_8$ of order $8$ is not an Oort group, by Example~2.12(f); and $Q_8$ is a subgroup of each semi-dihedral group or generalized quaternion group (e.g.\ by [As], p.115, Ex.~3(6)).  So $G$ cannot be a semi-dihedral group or generalized quaternion group, and the result follows. \qed

\medskip

Note that the Klein four group $C_2^2$, which is an Oort group and a local Oort group (see Example~2.12(e)), is included in Theorems~4.4 and 4.5 as the dihedral group $D_4$.  

\medskip

\noindent{\it Remark.}  In odd characteristic, Theorem~3.3 and Corollary~3.4 give the same necessary condition for being an Oort group or a local Oort group.  But in characteristic $2$, the necessary  condition in Theorem~4.4 to be a local Oort group is weaker than the corresponding condition to be an Oort group in Theorem~4.5.  These results suggest the question of whether, at least in odd characteristic, a cyclic-by-$p$ group is an Oort group if and only if it is a local Oort group.  The forward direction was shown in Corollary~2.6.  The converse is open, but it would follow in odd characteristic from Conjecture~1.1.  Namely, by that conjecture 
and Theorem~3.3, we need only consider local Oort groups $D_{2p^n}$.  By Proposition~2.11, $D_{2p^m}$ is a local Oort group for all $m \le n$.  By Conjecture~1.1 and Corollary~2.6, every cyclic group is a local Oort group.  So by Theorem~2.4, $D_{2p^n}$ is an Oort group, proving the converse for $p$ odd, assuming Conjecture~1.1.

\medskip

Applying Corollary~2.8 to Theorem~4.5 we obtain

\medskip

\noindent{\bf Corollary 4.6.} 
{\sl If $G$ is an Oort group in characteristic $2$, then every cyclic-by-$2$ subgroup of $G$ is isomorphic to a cyclic group, a dihedral group, or $A_4$.}

\medskip

Consequences of this result will be explored in [CGH2].

\medskip

\noindent{\bf Corollary 4.7.} 
{\sl Let $G$ be an Oort group for $k$, where $k$ has characteristic $2$.  Then the Sylow $2$-subgroups of $G$ are cyclic or dihedral.}

\medskip

\noindent{\it Proof.}  We proceed as in the proof of Corollary~3.6.  By Corollary~2.7, a Sylow $2$-subgroup $P \subset G$ is an Oort group for $k$.  Since $P$ is a $2$-group, it is not isomorphic to $A_4$.  So Corollary~4.5 implies that $P$ is cyclic or dihedral.  \qed

\medskip

\noindent{\bf Corollary 4.8.} 
{\sl Let $G$ be an Oort group for $k$, where  ${\rm char}\,k = 2$.  Let $P \subset G$ be a $2$-subgroup of $G$, with Frattini subgroup $\Phi$.  Suppose that $g \in G$ is an element of odd order that normalizes $P$ but does not centralize $P$.  Then $g$ has order $3$; $P$ has rank $2$; and the conjugation action of $g$ generates the automorphism group of $P/\Phi \approx C_2^2$.}

\medskip

\noindent{\it Proof.} Let $C$ be the subgroup generated by $g$.  
Since $g$ normalizes $P$, the subgroup $I=P.C$ generated by $P$ and $g$ is cyclic-by-$p$, and is an Oort group for $k$ by Corollary~2.7.  Since $g$ does not centralize $P$, it is not the identity element; and so $I$ strictly contains $P$ and is not a $2$-group.  Similarly, $I$ is not abelian.  So by Theorem~4.5, $I$ is isomorphic to $A_4$, and the conclusion follows.  \qed

\bigskip

\noindent{\bf References.}

[As] M.~Aschbacher.  ``Finite Group Theory'', second edition.  Cambridge Univ.\ Press, 2000.

[B] J.~ Bertin.  
 Obstructions locales au rel\`evement de rev\^etements galoisiens de courbes lisses.  C. R.  Acad.  Sci.  Paris S\'er. I Math.\  {\bf 326} (1998), no. 1, 55-58.
 
[BM] J.~Bertin, A.~M\'ezard. D\'eformations formelles de rev\^etements: un principe local-global.   Israel J. Math.  {\bf 155}  (2006), 281-307.

[BW] I.~Bouw, S.~Wewers.  The local lifting problem for dihedral groups. Duke Math.\ J.\ {\bf 134} (2006), 421-452.

[CGH1] T.~Chinburg, R.~Guralnick, D.~Harbater.  Bertin groups and local lifting problems.  To appear.

[CGH2] T.~Chinburg, R.~Guralnick, D.~Harbater.  On the structure of global Oort groups.  To appear.

[dJ] A.J.~de Jong. Families of curves and alterations.  Ann.\ Inst.\ Fourier (Grenoble) {\bf 47} (1997), 599-621.

[Go] D.~Gorenstein. Finite Groups, Harper  and Row, New York, 1968.

[GM] B.~Green, M.~Matignon. Liftings of Galois covers of smooth curves. Compositio Math., {\bf 113} (1998), 237-272.

[Grn] B.~Green. Automorphisms of formal power series rings over a valuation ring.  In: ``Valuation theory and its applications, Vol. II'', Fields Inst.\ Commun., vol.~33, AMS, 2003, pp. 79-87.

[Gr] A.~Grothendieck. ``Rev\^etements \'etales et groupe
fondamental'' (SGA 1). Lecture Notes in Mathematics, vol.\ 224,
Springer-Verlag, 1971.

[GS] R.~Guralnick, K.~Stevenson.  Prescribing ramification. 
In: ``Arithmetic fundamental groups and noncommutative algebra'' (M.~Fried, Y.~Ihara, eds.), AMS Proc. Symp. Pure Math.\ series, vol.~70, 2002, pp.~387-406.

[Ha1] D.~Harbater. Fundamental groups and embedding problems in characteristic $p$.  In ``Recent developments in the inverse
Galois problem'' (M.~Fried, et al., eds.), AMS Contemp.\ Math.\ Series, vol.~186, 1995, pp.~353-369.

[Ha2] D.~Harbater.  Embedding problems with local conditions.  Israel
J.\ of Math., {\bf 118} (2000), 317-355.

[HS] D.~Harbater, K.~Stevenson.  Patching and thickening problems.  J.\ Alg.\ {\bf 212} (1999), 272-304.

[H] R.~Hartshorne. ``Algebraic Geometry''.  Springer Graduate Texts in Mathematics, vol.~52, 1977.

[Ka] N.~Katz. Local-to-global extensions of representations of fundamental groups.  Ann.\ Inst.\ Fourier, Grenoble {\bf 36} (1986), 69-106.

[Ma] M.~Matignon.  $p$-groupes ab\'eliens de type $(p,...,p)$  et disques ouverts $p$-adiques.  Manuscripta Math.\ {\bf 99} (1999),  93-109.

[Oo] F.~Oort. Lifting algebraic curves, abelian varieties, and their endomorphisms to characteristic zero.  Proc.\ Symp.\ Pure Math., vol.~46, 1987.

[OSS] F.~Oort, T.~Sekiguchi, N.~Suwa. On the deformation of Artin-Schreier to Kummer.  Ann.\ Sci.\ \'Ecole Norm.\ Sup.\ {\bf 22} (1989), 345-375.

[Pa] G.~Pagot. Rel\`evement en caract\'eristique z\'ero d'actions de groupes ab\'eliens de type $(p,\dots,p)$. Ph.D.\ thesis, Universit\'e Bordeaux 1, 2002. 

[Po] F.~Pop. \'Etale Galois covers of affine smooth curves.    Invent.\ Math., {\bf 120} (1995), 555-578.

[Ra] M.~Raynaud.  $p$-groupes et r\'eduction semi-stable des courbes.  ``The Grothendieck Festschrift'', Vol. III, pp.~179-197,
Progr.\ Math., vol.~88, Birkh\"auser, Boston, 1990.

[Se] J.-P.~Serre.  ``Local Fields''.  Graduate Texts in Math., vol.~67, Springer-Verlag, 1979.

[Si] J.~Silverman.  ``The Arithmetic of Elliptic Curves''.  Graduate Texts in Math., vol.~106, Springer-Verlag, 1986.

[Su] M.~Suzuki. ``Group Theory II''. Grundlehren Math.\ series, vol.~248, Springer-Verlag, 1982.

[TT] O.~Taussky.
A remark on the class field tower. J.\ London Math.\ Soc.\ {\bf 12}, (1937), 82-85.

\medskip

\small 

\noindent T.~Chinburg: Department of Mathematics, University of Pennsylvania, Philadelphia, PA 19104-6395, USA; {\smtt ted@math.upenn.edu}

\noindent R.~Guralnick: Department of Mathematics, University of Southern California, Los Angeles, CA 90089-2532, USA; {\smtt guralnic@usc.edu}

\noindent D.~Harbater: Department of Mathematics, University of Pennsylvania, Philadelphia, PA 19104-6395, USA; {\smtt harbater@math.upenn.edu}

\end